\documentclass[11pt]{elsarticle}
\usepackage{amsmath,amssymb,amsthm}
\usepackage{mathtools, nccmath}
\usepackage{enumitem}
\usepackage[left=1.00in,right=1.00in,top=1.00in,bottom=1.00in]{geometry}

\usepackage[normalem]{ulem} 

\usepackage{graphicx}
\usepackage{mathtools}
 \usepackage{color} 

\newlength{\bibitemsep}
\newlength{\bibparskip}
\let\oldthebibliography\thebibliography
\renewcommand\thebibliography[1]{%
  \oldthebibliography{#1}%
  \setlength{\parskip}{\bibitemsep}%
  \setlength{\itemsep}{\bibparskip}%
}

\usepackage{titlesec}
\titlespacing*{\section}{0pt}{0.5\baselineskip}{0.3\baselineskip}

\journal{MS417 - Student Competition - USNCCM17}

\begin{document}

\begin{frontmatter}

\title{A Scalable Algorithm for Multi-Material Design of Thermal Insulation Components Under Uncertainty}

\author[A1]{Jingye Tan\corref{corauth}}
\cortext[corauth]{Student Contributor}
\ead{jtan32@buffalo.edu}


\author[A1]{Danial Faghihi}
\ead{danialfa@buffalo.edu}

\address[A1]{Department of Mechanical and Aerospace Engineering, 
University at Buffalo, Buffalo, NY, USA}

\begin{abstract}
This work presents a scalable computational framework for optimal design under uncertainty with application to multi-material insulation components of building envelopes.
The forward model consists of a multi-phase thermo-mechanical model of porous materials governed by coupled partial differential equations (PDEs). The design parameter (material porosity) is an uncertain and space-dependent field,
resulting in a high-dimensional PDE-constrained optimization under uncertainty problem after finite element discretization.
The robust design framework uses a risk-averse formulation consisting of the mean and variance of the design objective to achieve target thermal and mechanical performances and mitigate uncertainty. 
To ensure the efficiency and scalability of the solution, a second-order Taylor approximation of the mean and variance and the low-rank structure of the preconditioned Hessian of the design objective are leveraged, which uncovers the low effective dimension of the high-dimensional uncertain parameter space. 
Moreover, a gradient-based optimization method is implemented using the Lagrangian formalism to derive expressions for the gradient and Hessian with respect to the design and uncertain parameters. 
Finally, approximated $\ell_0$ regularization functions are utilized via a continuation numerical scheme to promote sparsity in the designed porosity. The framework's accuracy, efficiency, and scalability are demonstrated with numerical examples of a building envelope insulation scenario.
\end{abstract}


\end{frontmatter}

\section{Introduction}
\label{sect:introduction}

The simulation-based design of multi-material components governed by partial differential equations (PDEs) models is challenging due to the high dimensionality of the design and uncertain parameters, the need to control the interface thickness between materials, and the costly solution of complex PDEs in large-scale systems. 
Moreover, the inevitable uncertainties in the model parameters (due to experimental noise and model inadequacy) and the design parameters (introduced during material fabrication) lead to computationally prohibitive optimization problems.


%

This work develops a scalable framework for PDE-constrained optimal multi-material design under uncertainty (DUU). The framework is applied to insulation components of building envelopes to achieve desired thermal insulation and mechanical strength.
The rest of the paper is organized as follows: The thermo-mechanical model is presented in Section 2.
The main components of the proposed DUU framework are summarized in Section 3. 
The framework's efficiency and scalability are demonstrated in Section 4 via numerical experiments. 
Section 5 presents conclusions and our ongoing efforts to extend the method and numerical experiments.

\section{Multiphase Thermomechanical Model of Porous Materials}


    We aim to design insulation components (thermal breaks) in building envelopes made of silica aerogel, a meso-porous material with super-insulation property.
    The thermo-mechanical behavior of silica aerogel is modeled using a continuum mixture theory \cite{tan2022predictive}. The governing equations under steady state thermal transport and static conditions reduce to the following PDEs,
            \vspace{-0.05in}
            \begin{equation}
                \begin{cases}
                    -\nabla \cdot(\phi_s \, \kappa_s \, \nabla T_s) &= -h(T_s-T_f),
                    \\
                    -\nabla \cdot(\phi_f \, \kappa_f \, \nabla T_f) &= h(T_s-T_f),
                \end{cases}
                \quad\quad
                ,
                \quad\quad
                \begin{cases}
                    C\,p = -(\nabla\cdot\mathbf{u}_s)  ,
                    \\
                    \nabla\cdot\mathbf{T}^\prime_s  +  (2\,\phi_f - 1)\,\nabla p  = 0   ,
                \end{cases}
                \label{pdes}
                \vspace{-0.05in}
            \end{equation}
            where, the state variables $T_s$, $T_f$, $\mathbf{u}_s$ and $p$ are solid and fluid temperatures, solid displacement and fluid pressure, respectively. The effective solid stress is $\mathbf{T}^\prime_s = 2\,\mu\,\boldsymbol{\epsilon}_s + \lambda\,\mathrm{tr}(\boldsymbol{\epsilon}_s)\,\mathbf{I}$ with $\boldsymbol{\epsilon}_s = \frac{1}{2}\big(  \nabla\mathbf{u}_s + (\nabla\mathbf{u}_s)^T  \big)$
            being the solid strain.
            The model parameters $\boldsymbol{\theta} =($$\kappa_s$, $\kappa_f$, $h$, $K$, $C$,  $\mu)$ can be determined from experimental data \cite{tan2022predictive}. The spatial distribution of the fluid volume fraction (porosity) $\phi_f = 1 - \phi_s$ is the design variable.
            The quantities of interests (QoIs) representing insulation $Q_T$ and mechanical $Q_M$ performances, with $j$ convection and $k$ traction boundaries, are
            \vspace{-0.05in}
            \begin{equation}
                \begin{cases}
                    Q_T 
                    &= \  
                    \frac{1}{2}\sum_{i = s,f} \big<\phi_i\,\kappa_i\,\nabla T_i, \nabla T_i \big>_{\Omega} 
                    + 
                    \sum_{i = s,f} \sum_{j} \big< \phi_i ( \frac{1}{2} T_i + T^\mathrm{amb}_j)   ,  T_i   \big>_{\Gamma_j},
                    \\
                    Q_M 
                    &= \ 
                    \frac{1}{2} \big<\mathbf{T}^\prime_s,\boldsymbol{\epsilon}_s \big> + \sum_{k} \big< \mathbf{t}, \mathbf{u}_s \big>_{\Gamma_k},
                \end{cases}
                \label{qois}
                \vspace{-0.05in} 
            \end{equation}
            %
            where $\big<\cdot,\cdot\big>$ represents integration of the two terms over the domain $\Omega$ or on the boundary $\Gamma$.

\section{A Scalable Framewrok for Optimal Design Under Uncertainty (DUU)}

    Despite their excellent insulation properties, silica aerogel's high porosity results in mechanical instability.
    We aim to determine the spatial distribution of porosity within the thermal break to achieve desired insulation and mechanical properties. 
    Nevertheless, uncertainty is associated with the porosity value and its spatial distribution.
    To model the uncertainty, 
    we consider the design parameter $d(\mathbf{x}) \in [0, 1]$ to be a linear map $f(\cdot)$ interpolates to silica aerogel porosity ($\phi_f(\mathbf{x})=[0.1, 0.9]$), 
    including uncertainty, i.e., $\phi_f \equiv f(d+m)$. 
    The uncertain parameter $m(\mathbf{x})$ is represented as Gaussian random field with the Mat\'{e}rn Coveriance kernel $\mathcal{C}=\mathcal{A}^{-2}$ \cite{lindgren2011explicit}, such that 

    \vspace{-0.05in}
    \begin{equation}
        \mathcal{A} \, m = \begin{cases}
            \gamma \, \nabla \cdot( \boldsymbol{\Theta} \,\nabla m ) + \delta \, m  \quad\mathrm{in\ }\Omega
                    \\
            (\boldsymbol{\Theta}\,\nabla m) \cdot \mathbf{n}  + \frac{\sqrt{\delta\,\gamma}}{1.42}\,m  \quad\mathrm{on\ }\Gamma
        \end{cases}
        ,
        \vspace{-0.05in}
    \end{equation}
    \useshortskip
    where $\gamma$ and $\delta$ control the variance and spatial correlation, and $\boldsymbol{\Theta}$ is an anisotropic tensor.



            Since the thermal and mechanical QoIs (design objective) depends on uncertain parameter $m$, we employ       
            a risk-averse PDE-constrained optimization formulation \cite{alexanderian2017mean,chen2019taylor} for the DUU problem,


            \vspace{0.05in}
            \begin{equation}
                    \min_d  \mathcal{J}(d) = \mathbb{E}[Q(d,m)] + \beta_V \,\mathbb{V}[Q(d,m)] + \beta_R \, R(d)
                    , \quad
                    \mathrm{subjected\ to}\mathrm{\  \eqref{pdes} },\  d \in [0,1],
                \label{cost}
                \vspace{-0.05in}
            \end{equation}

\noindent
            where $\mathbb{E}[\cdot]$ and $\mathbb{V}[\cdot]$ represent mean and variance measures, and $R(d)$ is a regularization function.  
            The design objective is $Q = \beta_M\,Q_M - Q_T$, where $\beta_M$ is the weight of the mechanical QoI relative to that of thermal, and 
            $\beta_V$ and $\beta_R$ are weights of the variance and regularization terms.

    \subsection{Taylor approximation of design objective}

        The Monte Carlo (MC) estimations of the mean and variance in \eqref{cost} involves solving the governing PDEs \eqref{pdes} for a high-volume samples of $m$ for each optimization iterations.
        To alleviate such prohibitive cost, 
        we develop an efficient algorithm by employing a second order Taylor approximation of the mean and variance \cite{alexanderian2017mean,chen2019taylor},
        evaluated at the uncertain parameter mean $\bar{m}$,

  \useshortskip
        \begin{equation}
            \mathbb{E}[Q] \approx \Bar{Q} + \frac{1}{2}
            \mathrm{tr}\big(\Bar{\mathcal{H}}_c\big)
            ,\,
            \mathbb{V}[Q] \approx \big<\Bar{Q}_m, \mathcal{C}\,\Bar{Q}_m \big> + \frac{1}{2}
            \mathrm{tr}\big(\Bar{\mathcal{H}}_c^2\big),
            \;
            \mathrm{with} \;\;
            \mathrm{tr}\big(\Bar{\mathcal{H}}_c\big) \approx \sum_{n\geq1}^N  \lambda_n
            ,\,
            \mathrm{tr}\big(\Bar{\mathcal{H}}_c^2\big) \approx \sum_{n\geq1}^N  \lambda_n^2.
            \label{meanvar}
            \end{equation}
        Above, $\Bar{Q}$, $\Bar{Q}_m$, $\Bar{\mathcal{H}}_c = \mathcal{C}\,\bar{\mathcal{H}}$ are the design objective, its gradient, and its covariance preconditioned Hessian, evaluated at $\bar{m}$.
        The trace terms $\mathrm{tr}(\cdot)$, are approximated with $N$ dominant eigenvalues $(\lambda_n)_{n\geq1}$ of the covariance-preconditioned Hessian $\bar{\mathcal{H}}_c$, due to their rapid decay (see Sec.\ref{sec4} Fig.\ref{result1}(b)),
        and $\lambda_n$ are obtain by solving the generalized eigenvalue problem
        using a randomized algorithm \cite{saibaba2016randomized}.
    
    \subsection{Gradient-based optimization}

        We employ the limited-memory Broyden–Fletcher–Goldfarb–Shanno optimization algorithm.
        We use Lagrangian formalism to derive the required derivatives, i.e.,
        $m$-gradient and $m$-Hessian of $Q$ in \eqref{meanvar} and $d$-gradient of the quadratic approximation of $\mathcal{J}(d)$ in \eqref{cost}.
        

    \subsection{Sparsity-enforcing regularization via continuation scheme}

            To control the interface thickness between two materials, we propose regularization term in \eqref{cost} as a weighted sum of the Tikhonov and approximated $\ell_0$ via a continuation scheme,  


            \vspace{-0.05in}
            \begin{equation}
            R(d)
            =
             \int_\Omega 
             \beta_{R_\mathrm{tik}} 
             \big(\nabla d : \nabla d \big)
             +
             \beta_{R_0}  
             \big(f_{\epsilon_0}(d)\big)
             \, \mathrm{d}\Omega, 
             \ 
             \mathrm{where,\ } f_{\epsilon_0}(d) = 
                \begin{cases}
                d/\epsilon_0,  &  d\in [  0,   \epsilon_0/2 ]
                \\
                p_{3}(d,\epsilon_0),  &  d\in ( \epsilon_0/2,  2\epsilon_0 )
                \\
                1, &  d\in    (2\epsilon_0,   1 ]
                \end{cases},
                \label{reg}
            \vspace{-0.05in}
            \end{equation}
            where $p_{3}(d,\epsilon_0)$ is a third order polynomial on $d$ uniquely determined 
            by $\epsilon_0 \in (0,1)$ to make $f_{\epsilon_0}(d)$ continuously differentiable. 
            %
            The continuation numerical scheme consists of
            first performing optimization with $\beta_{R_0}=0$ and $\beta_{R_\mathrm{tik}} =1$ for fast convergence, resulting in a possibly non-sparse optimal design $d_\mathrm{opt}^\mathrm{tik}$.  
            Second, we set $\beta_{R_0}=1$ and $\beta_{R_\mathrm{tik}} =0$, 
            and with a decreasing sequence of 
            $\big\{\epsilon_0^{(i)} = (\frac{1}{2})^i \big\}_{i=1}^{K}$ 
            for $K$ continuation iterations.
            In the first iteration, $d_\mathrm{opt}^\mathrm{tik}$ is used to initialize the continuation with $\epsilon_0^{(1)}$ to obtain $d_\mathrm{opt}^{(1)}$.
            Each subsequence optimization with $\epsilon_0^{(i+1)}$ will be starting from the previous optimal design $d_\mathrm{opt}^{(i)}$. 
            Within this scheme, the mesh can also be adaptively refined in non-sparse regions to 
            efficiently control materials interface thicknesses.

\section{Numerical Results} \label{sec4}

We present a set of preliminary results of the proposed DUU framework on a L-shape insulation component, exposed to convection with a hot temperature and uniform traction load on the outer boundary, and convection to a cold temperature and fixed displacement on the inner boundary.

\begin{figure}[h]
\vspace{-0.1in}
    \centering
    \begin{center}
    \includegraphics[width=1.0\textwidth, clip=true, trim = 0mm 0mm 0mm 0mm]{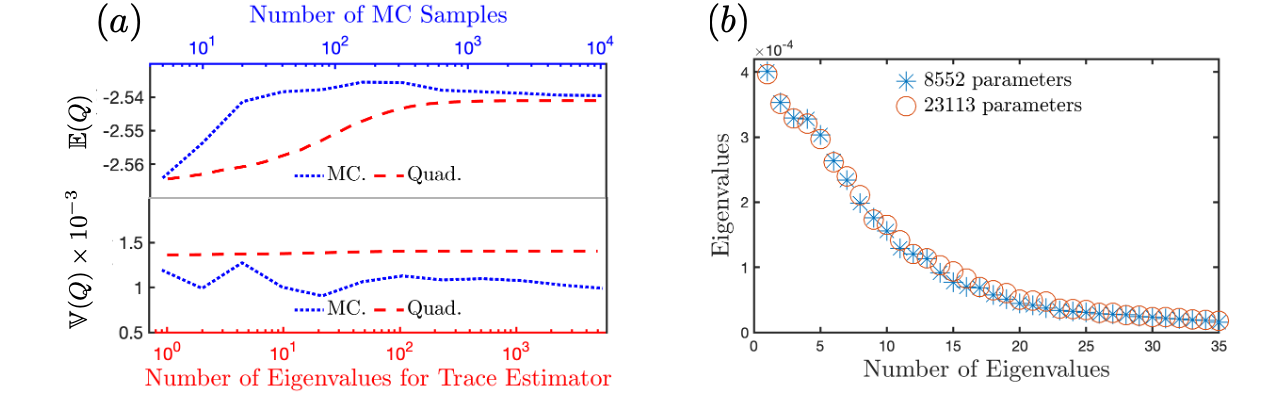}\\
    \end{center}
    \vspace{-0.25in}
    \caption{
    (a) Convergence of the quadratic Taylor approximation (Quad) of $Q$ with respect to the number of dominant eigenvalues in comparison with the Monte Carlo estimation (MC).
    (b) Decay of the eigenvalues of the covariance-preconditioned Hessian with different uncertain parameter dimensions (mesh discretizations).
%
}
    \label{result1}
\end{figure}

\begin{figure}[htpb]
\vspace{-0.1in}
    \centering
    \begin{minipage}{0.6\textwidth}
    \includegraphics[width=\textwidth, clip=true, trim = 0mm 0mm 0mm 0mm]{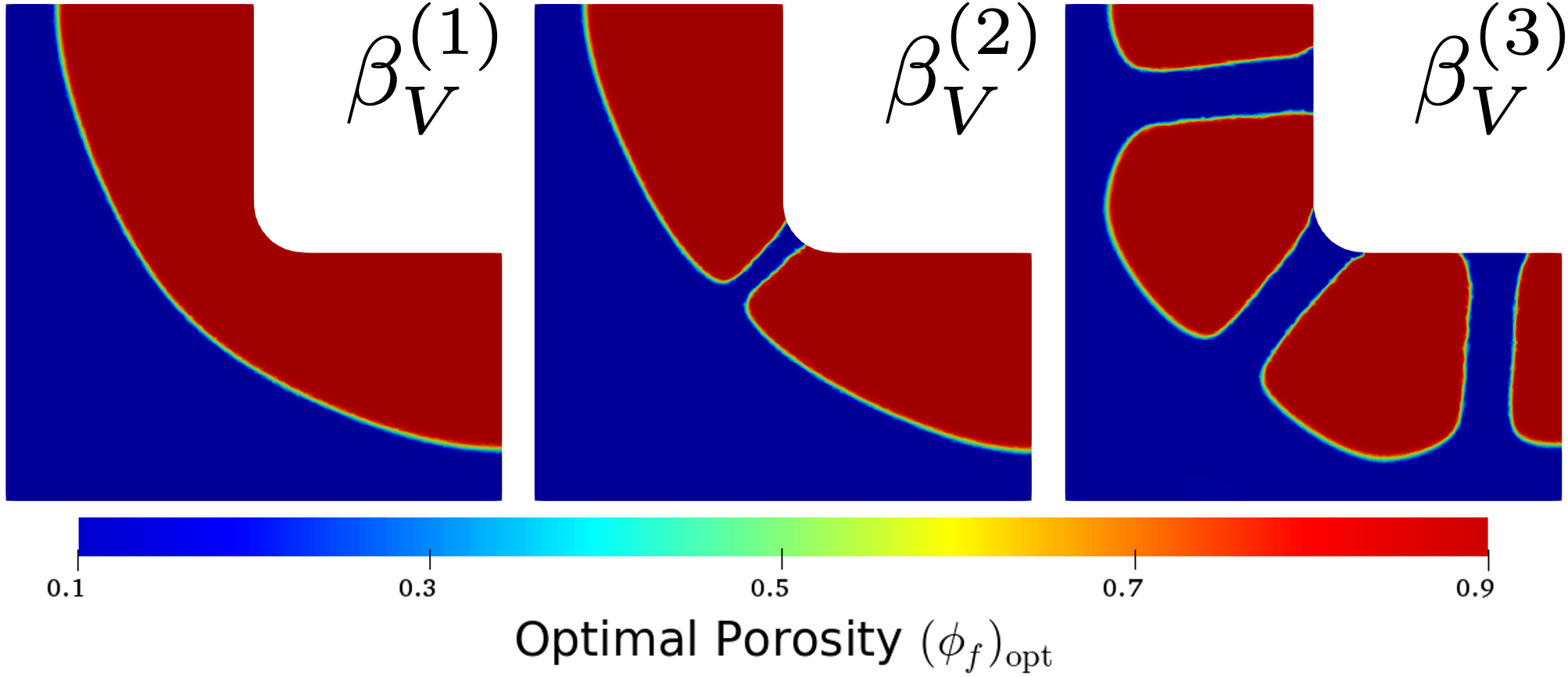}
    \end{minipage}
    \hfill
    \begin{minipage}{0.39\textwidth}
    \includegraphics[width=\textwidth, clip=true, trim = 0mm 0mm 0mm 0mm]{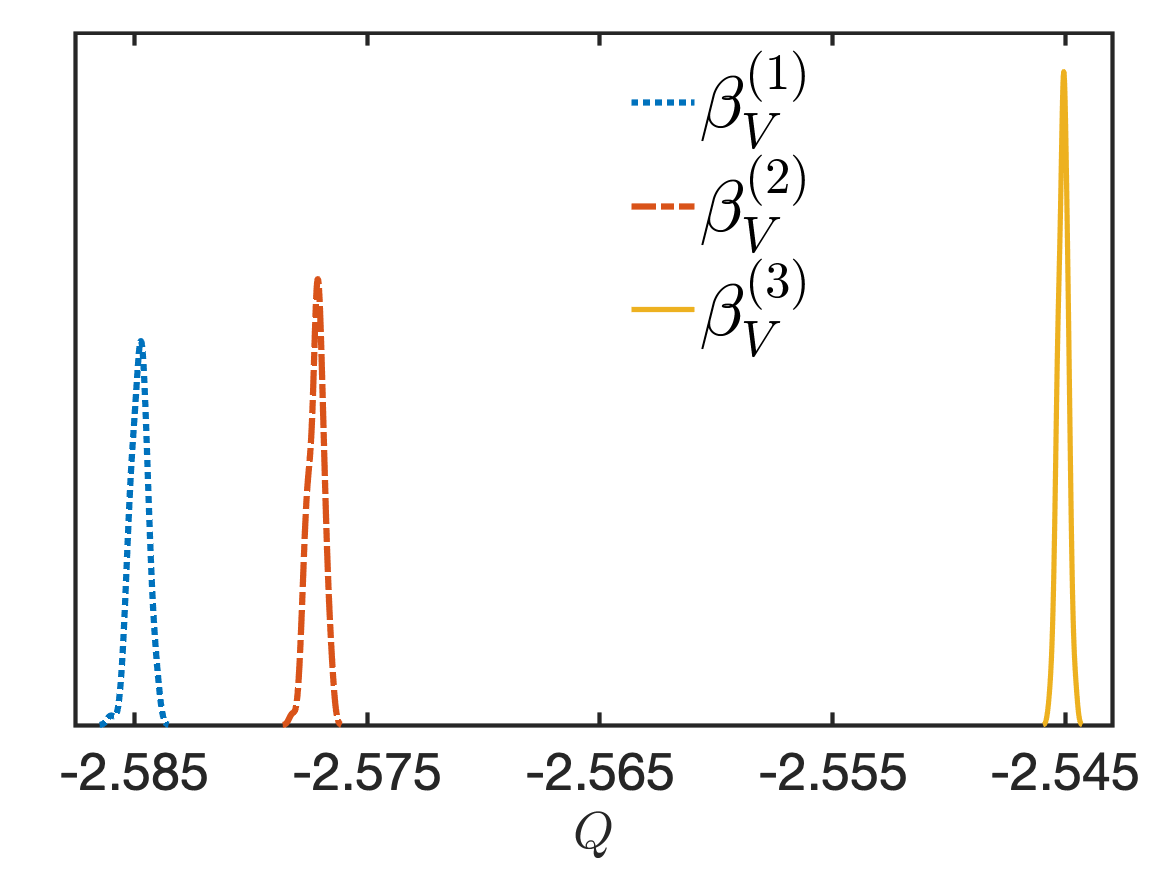}
    \end{minipage}
    \vspace{-0.2in}
    \caption{Optimal design for different variance weights, $\beta_V^{(1)}=0$, $\beta_V^{(2)}=1$E5, and $\beta_V^{(3)}=1$E6
    and the distributions of design objective $Q$. 
    The variances of $Q$ are $1.3\times 10^{-7}$, $1.1\times 10^{-7}$, and $4.7\times 10^{-8}$ for $\beta_V^{(1)}$, $\beta_V^{(2)}$, $\beta_V^{(3)}$, respectively.}
    \label{result2}
    \vspace{-0.15in}
\end{figure}

Fig \ref{result1}(a) shows the convergence of the quadratic approximation of $Q$ in \eqref{meanvar} with respect to the number of dominant eigenvalues $\lambda_n$. 
The results indicate that using the first 25 dominant eigenvalues, the quadratic approximation achieves below $0.5\%$ relative error compared to MC estimation (with 10240 samples).
Fig \ref{result1}(b) represents the fast decay of the eigenvalues for two finite element discretizations. This indicates that the  
computational work (required number of PDE solves) of the quadratic approximation does not depend on the high dimension of the uncertain parameter, but rather on the effective low-dimension stem from the rank of the preconditioned Hessian, providing the scalability of the algorithm.
The optimal design results are shown in Fig \ref{result2}
with three different levels of variance weights in \eqref{cost}.
As $\beta_V$ increases, the spatial patterns of optimal porosity provide more mechanical support, while the variance of $Q$ reduces by trading off the mean performance.

\section{Conclusion}
\label{sect:conclusion}

We have developed a scalable algorithm for optimally distributing porosity in aerogel-based insulation components. Our approach includes a risk-averse formulation for robust design and a Taylor approximation of the mean and variance of the design objective.
We utilize a randomized algorithm to estimate the trace of the preconditioned Hessian of the design objective and a gradient-based optimization method using Lagrangian formalism.
Our preliminary numerical experiments show that the framework's computational work (required number of PDE solves) does not depend on the discretized high-dimensional uncertain parameter, but rather on the effective low-dimension stem from the rank of the preconditioned Hessian, providing the scalability of the algorithm.

Our ongoing work complements the results presented in this paper in two ways. 
First, we use the quadratic approximation of the design objective as a control variate, leading to significant computational savings compared to a plain Monte Carlo method and increased numerical accuracy compared to quadratic approximation alone.
Second, we apply the framework to a more complex scenario involving aerogel thermal breaks in building envelopes and column systems, requiring a seamless blend of the continuation scheme for $\ell_0$-sparsified designs with adaptive mesh refinement to control the interface thickness between two different material porosity.

\section{Acknowledgments}

\small
The support from the U.S. National Science Foundation CAREER Award CMMI-2143662 is acknowledged.


\small

\bibliographystyle{elsarticle-num} 
\bibliography{biblio.bib}

\normalsize








\end{document}